\documentclass[12pt]{article}
\usepackage{amsfonts, amsmath, amstext, amssymb}

\usepackage[dvips]{graphicx}
\usepackage{latexsym}

\newtheorem{thm}{Theorem}[section]
\newtheorem{prop}[thm]{Proposition}

\newtheorem{lemma}[thm]{Lemma}

\newtheorem{cor}[thm]{Corollary}

\newtheorem{definitiontemp}[thm]{Definition}
\newenvironment{defn}{\begin{definitiontemp}
\normalfont}{\end{definitiontemp}}

\def\bec{\begin{cor}}
\def\enc{\end{cor}}

\def\ep{\varepsilon}
\newcommand{\HP}{\emptyset'}

\def\bet{\begin{thm}}
\def\ent{\end{thm}}

\def\becor{\begin{cor}}
\def\encor{\end{cor}}

\def\bel{\begin{lem}}
\def\enl{\end{lem}}

\def\bedef{\begin{defn}}
\def\endef{\end{defn}}

\def\bep{\begin{prop}}
\def\enp{\end{prop}}

\newenvironment{pf}{\begin{trivlist}\item[\hskip\labelsep
{\it Proof.}]}{\end{trivlist}}

\newenvironment{pftitle}[1]{\begin{trivlist}\item[\hskip\labelsep
{\it #1.}]}{\end{trivlist}}

\newcommand{\set}[2]{\ensuremath{ \{ #1 : #2 \} }}

\newcommand{\N}{\mathbb{N}}
\newcommand{\Z}{\mathbb{Z}}

\newcommand{\Q}{\mathbb{Q}}
\newcommand{\R}{\mathcal{R}}
\newcommand{\C}{\mathcal{C}}
\newcommand{\A}{\mathcal{A}}
\newcommand{\B}{\mathcal{B}}

\renewcommand{\P}{\mathbb{P}}
\newcommand{\U}{\mathcal{U}}
\newcommand{\V}{\mathcal{V}}
\newcommand{\W}{\mathcal{W}}

\newcommand{\T}{\mathcal{T}}

\newcommand{\qvec}{\vec{q}}
\newcommand{\rvec}{\vec{r}}
\newcommand{\svec}{\vec{s}}
\newcommand{\tvec}{\vec{t}}

\newcommand{\wvec}{\vec{w}}

\newcommand{\Tvec}{\vec{T}}

\newcommand{\Xvec}{\vec{X}}
\newcommand{\yvec}{\vec{y}}
\newcommand{\Yvec}{\vec{Y}}
\newcommand{\zvec}{\vec{z}}
\newcommand{\Zvec}{\vec{Z}}

\newcommand{\initial}{\sqsubseteq}

\newcommand{\la}{\langle}
\newcommand{\ra}{\rangle}

\newcommand{\HTP}[1]{HTP(#1)}
\newcommand{\HTPQ}{HTP(\Q)}
\newcommand{\HTPZ}{HTP(\Z)}
\newcommand{\ldotc}{,\ldots,}

\def\diverges{\!\uparrow}
\def\converges{\!\downarrow}
\newcommand{\at}{\char'100}

\newcommand{\Leg}[2]{\left(\begin{smallmatrix} \underline{#1} \\ #2 \end{smallmatrix}\right)}

\newcommand{\qed}{\hbox to 0pt{}\nobreak\hfill\rule{2mm}{2mm}}

\newcommand{\dom}[1]{\text{dom}(#1)}

\newcommand{\bfd}{\boldsymbol{d}}
\def\bfz{\boldsymbol{0}}
\def\s01{\ensuremath{\Sigma^0_1}}
\def\d02{\ensuremath{\Delta^0_2}}
\def\phi{\varphi}

\def\res{\!\!\upharpoonright\!}

\newcommand{\comment}[1]{}

\begin{document}

\title{HTP-Complete Rings of Rational Numbers}
\author{Russell Miller
\thanks{The author was partially supported by Grant \# 581896
from the Simons Foundation, and by several grants from
the PSC-CUNY Research Award Program.  He offers sincere thanks
to the anonymous referee for several important corrections and clarifications.
}
}
\maketitle

\begin{abstract}
For a ring $R$, Hilbert's Tenth Problem $HTP(R)$ is the set of polynomial
equations over $R$, in several variables, with solutions in $R$.
We view $HTP$ as an enumeration operator, mapping each set $W$ of prime numbers
to $HTP(\Z[W^{-1}])$, which is naturally viewed as a set of polynomials
in $\Z[X_1,X_2,\ldots]$.  It is known that for almost all $W$, the jump $W'$
does not $1$-reduce to $HTP(R_W)$.  In contrast, we show that every Turing degree
contains a set $W$ for which such a $1$-reduction does hold:
these $W$ are said to be \emph{HTP-complete}.
Continuing, we derive additional results regarding the impossibility
that a decision procedure for $W'$ from $HTP(\Z[W^{-1}])$ can succeed
uniformly on a set of measure $1$, and regarding the consequences
for the boundary sets of the $HTP$ operator in case $\Z$ has
an existential definition in $\Q$.
\end{abstract}

\noindent
$\textbf{Key Words}$: boundary sets, computability theory, Hilbert's Tenth Problem, HTP-completeness, subrings of the rational numbers.

\section{Introduction}
\label{sec:intro}

For a ring $R$, Hilbert's Tenth Problem $HTP(R)$ is the set of polynomial
equations over $R$, in several variables, with solutions in $R$:
$$ HTP(R)=\bigcup_{n\in\N}\set{f\in R[X_1,\ldots,X_n]}{(\exists x_1\ldotc x_n\in R)~f(x_1\ldotc x_n)=0}.$$
For countable rings $R$, one can ask effectiveness questions about
$HTP(R)$, which is always computably enumerable relative to a
presentation of $R$ (that is, relative to the atomic diagram of
a ring isomorphic to $R$ with $\omega$ as its underlying set; cf.\ \cite{M08}).
That it may be \emph{properly} computably enumerable was established
for the fundamental example $R=\Z$ by Matiyasevich in \cite{M70},
completing work by Davis, Putnam, and Robinson in \cite{DPR61}:
they showed that the Halting Problem is $1$-reducible to $HTP(\Z)$.
This was the resolution of Hilbert's original Tenth Problem, in which
Hilbert had demanded an algorithm deciding membership in $HTP(\Z)$.
(Simpler constructions exist of computable rings $R$ for which
$HTP(R)$ is undecidable; indeed one will be described at the end
of Section \ref{sec:construction}.)  In contrast, the decidability
of $HTP(\Q)$ remains an open question, and is the subject of intense study.

The subrings of $\Q$
are in bijection with the subsets $W$ of the set $\P$ of all prime numbers,
via the map $W\mapsto\Z[W^{-1}]$.  Thus these subrings form a topological
space homeomorphic to Cantor space $2^{\P}$, on which one
can therefore consider questions of Lebesgue measure and Baire category.
We then view $HTP$ as an operator, mapping each $W\subseteq\P$ to the
set $HTP(\Z[W^{-1}])$.  Notice that, to decide $HTP(\Z[W^{-1}])$, one need
only decide membership in it for polynomials from $\Z[X_1,X_2,\ldots]$, and
so we normally view $HTP(\Z[W^{-1}])$ as its own intersection with
$\Z[X_1,X_2,\ldots]$.  This allows for a uniform G\"odel coding of each
$HTP(\Z[W^{-1}])$ as a subset of $\omega$.  For simplicity we often
write $R_W$ for the ring $\Z[W^{-1}]$, and $HTP(R_W)$ for its $HTP$.

This article continues a program by the author of approaching $HTP(\Q)$
by viewing the collection of all subrings of $\Q$ as a topological space
in this way and considering the ``common'' behavior of the sets $HTP(R)$.
Certain properties hold of $HTP(R)$ for every $R$ in a ``large'' set of rings
(corresponding to a subset of measure $1$ within Cantor space, say,
or to a comeager subset), while other properties occur less frequently.
In turn, that program fits within the broader framework of examining $HTP(R)$
for subrings $R\subseteq\Q$ in general, an endeavor that includes the notable work
of Poonen \cite{P03,P09} and many others.  The results in Section
\ref{sec:construction} here may be seen as fitting into the larger program,
whereas the results in Sections \ref{sec:enumops} and
\ref{sec:definability} are mainly of interest within the smaller program.

The jump $W'$ of a set $W\subseteq\P$ is readily seen to give an upper bound
for the complexity of the set $HTP(R_W)$, which is always c.e.\ relative
to $W$.  It is known from work in \cite{KM19}
that for many subrings $R_W$ of $\Q$, the complexity of $HTP(R_W)$
is strictly below that of $W'$; this will be generalized here in Theorem \ref{thm:JK}.
In contrast, in Section \ref{sec:construction}, we will show that sets $W$
with $HTP(R_W)\equiv_1 W'$ are ubiquitous, in the sense that they exist
within every Turing degree.  Such sets will be said to be \emph{HTP-complete}.
In Section \ref{sec:enumops} we will consider the weaker relationship
of Turing-reducibility between $HTP(R_W)$ and $W'$, whereas in Section
\ref{sec:definability} we will consider consequences that would follow in case
$\Z$ has an existential definition in $\Q$ -- which is an open question, and
represents the strongest conjecture normally considered about the
difficulty of deciding $HTP(\Q)$.  (If such an existential definition exists,
then $HTP(\Z)$ itself $1$-reduces to $HTP(\Q)$, and therefore so does
the Halting Problem.)  Before that, Sections \ref{sec:background} and
\ref{sec:HTPcomplete} provide lemmas established earlier
in \cite{EMPS15} and \cite{KM19}, which will be of use
in the subsequent sections.  For basic information about computability theory,
\cite{S87} remains an excellent source, but \cite{R67} will be more helpful
for the concept of enumeration reducibility in Section \ref{sec:enumops}.

\section{Background}
\label{sec:background}

The following lemma appears in \cite{EMPS15}, but it has been known
ever since the pioneering work of Julia Robinson in \cite{R49}, in which
Robinson gave the first definition of the set $\Z$ within the field $\Q$.
It will be important for us at several junctures, as it enables us to ``ignore''
a given finite set of primes when dealing with the $HTP$ operator, and thus
facilitates finite-injury constructions.

\begin{prop}[Robinson \cite{R49}]
\label{prop:semilocal}
For every prime $p$, there is
a polynomial $g_{p}(Z,X_1,X_2,X_3)$ such that for all rationals $q$, we have
$$ q\in R_{(\P-\{p\})} \iff g_{p}(q,\Xvec)\in\HTPQ.$$
Moreover, $g_p$ may be computed uniformly in $p$.
\end{prop}
\begin{cor}
\label{cor:semilocal}
For each finite subset $A_0\subseteq\P$, a polynomial
$f(Z_0,\ldots,Z_n)$ has a solution in $R_{(\P-A_0)}$ if and only if
$$(f(\Zvec))^2 + \sum_{p\in A_0, j\leq n} (g_{p}(Z_j,X_{1,j,p},X_{2,j,p},X_{3,j,p}))^2$$
has a solution in $\Q$.
\end{cor}
For a proof using the more recent results of \cite{K15},
see Proposition 5.4 in \cite{EMPS15}.

Our principal tool for proving Theorem \ref{thm:construction} and its corollaries
will be the equations $X^2+qY^2=1$.  Here we review the relevant
number-theoretic results, previously applied to this purpose in \cite{KM19} by Kramer
and the author.  The main point is that for each odd prime $q$, there is an
infinite decidable set $V$ of primes such that $R_V$ contains no nontrivial solutions
to $X^2+qY^2=1$, yet for every $p\notin V$, the ring $\Z[\frac1p]$ does contain a
nontrivial solution.  (Here the trivial solutions are $(\pm 1,0)$, which in Section \ref{sec:construction}
will be ruled out as solutions, at the cost of turning $(X^2+qY^2-1)$ into a messier polynomial.)

\begin{defn}
\label{defn:q-appropriate}
For a fixed odd prime $q$, a prime $p$ is \emph{$q$-appropriate}
if $p$ is odd and $p\neq q$ and $\Leg{-q}p=1$ (that is, $-q$ is a square modulo $p$).
\end{defn}

The crux of Lemma \ref{lemma:qpolys} is that the $q$-appropriate primes
are precisely the possible factors of the denominator in a nontrivial solution
to $x^2+qy^2=1$, thus justifying the term \emph{$q$-appropriate}.
This lemma comprises Lemmas 4.2 and 4.4 from \cite{KM19},
where the proofs may be found.
\begin{lemma}
\label{lemma:qpolys}
Fix an odd prime $q$, and let $x$ and $y$ be positive rational numbers with $x^2+qy^2=1$.
Then every odd prime factor $p$ of the least common denominator $c$ of $x$ and $y$ must be $q$-appropriate.

Conversely, suppose that $p$ is $q$-appropriate.
Then there is a primitive solution $(a,b,p^k)$ to $X^2 + qY^2 = Z^2$
with $k \ge 1$.  Hence there is a nontrivial solution to $X^2+qY^2 = 1$
in the ring $\Z[\frac1p]$.

Finally, for $q\equiv 3\bmod 4$, a prime
$p\neq q$ is $q$-appropriate if and only if $p$ is a square modulo $q$.
When $q\equiv 1\bmod 4$, a prime $p\neq q$ is $q$-appropriate
if and only if one of the following holds:
\begin{itemize}
\item
$p\equiv 1\bmod 4$ and $p$ is a square modulo $q$.
\item
$p\equiv 3\bmod 4$ and $p$ is not a square modulo $q$.
\end{itemize}
It follows that $q$-appropriateness of $p$ is decidable uniformly in $q$.
\qed\end{lemma}

Of course, a single prime $p$ can be $q$-appropriate for many different $q$.
Therefore, adjoining $\frac1p$ to a ring may create solutions of the equations
$X^2+qY^2=1$ for many values of $q$.  The purpose of the following corollary
(in concert with Proposition \ref{prop:semilocal} above) is to allow us
to create a solution to the equation of our choice without disrupting
the solvability of other such equations in the ring.
Corollary \ref{cor:qpolys} is a modest extension of Corollary 4.3
in \cite{KM19}, where $I$ was the set $\{ 0,1,\ldots,e-1\}$.

\begin{cor}
\label{cor:qpolys}
Let $3=q_0<q_1<\cdots$ be the odd prime numbers.
Then, for every $e\in\omega$ and every finite set $I\subseteq\omega-\{ e\}$,
there are infinitely many primes $p$ that are $q_e$-appropriate but (for all $i\in I$)
are not $q_i$-appropriate.
\end{cor}
\begin{pf}
Write $I=\{i_0 < \cdots < i_j\}$.  Our goal is to show that there is a residue $n$
modulo $m=4q_{i_0}\cdots q_{i_j}\cdot q_e$ which is prime to $m$ and
satisfies all the criteria dictated by the final part of Lemma \ref{lemma:qpolys},
so that each prime $p$ congruent to $n$ modulo $m$ will satisfy the corollary.
For $m_0=4q_{i_0}$, we choose $n_0\equiv 1\bmod 4$ such that $n_0$
is not a square mod $q_{i_0}$, noting that for each residue $r$ mod $q_{i_0}$,
one of $r, r+q_{i_0}, r+2q_{i_0}, r+3q_{i_0}$ must be $\equiv 1\bmod 4$.
Setting $m_{k+1}=m_k\cdot q_{i_{k+1}}$ inductively for $k<j$, we see that
the residue $n_k$ mod $m_k$ yields distinct residues $n_k, n_k+m_k,
\ldots, n_k+m_k(q_{i_{k+1}}-1)$ mod $m_{k+1}$, each congruent to a distinct
residue mod $q_{i_{k+1}}$.  So we may choose $n_{k+1}$ to be one
of these residues which is not a square mod $q_{i_{k+1}}$.  Once we
have produced $n_j$ (a residue modulo $m_j$), we do the same process
with $q_e$, except that now we choose the new residue $n$ mod $m$
so that $n$ is a nonzero square mod $q_e$.  With $n\equiv 1\bmod 4$, this means
that each prime with residue $n$ mod $m$ will be $q_e$-appropriate
(as $n$ is a square mod $q_e$), but will be $q_{i_k}$-inappropriate
for all $k\leq j$ (as $n$ is not a square mod $q_{i_k}$).  Finally, $n$
is nonzero modulo $4$, modulo $q_e$, and modulo each $q_{i_k}$,
hence is prime to $m$.
Therefore, Dirichlet's theorem on arithmetic progressions
(see \cite[Chap.\ 6, \S 4]{S73}) shows that infinitely many primes
are congruent to $n$ mod $m$, so the corollary holds.
\qed\end{pf}

The equation $X^2+qY^2=1$ may be seen as stating that
an element $x+y\sqrt{-q}$ of the field $\Q(\sqrt{-q})$ has norm $1$ there.
It has been pointed out that many other sets of norm equations
of totally complex extensions of degree $2$ would admit
similar results to Lemma \ref{lemma:qpolys} and Corollary
\ref{cor:qpolys}, and therefore could presumably be used
to give alternative constructions for Theorem \ref{thm:construction}.

\section{HTP-Completeness}
\label{sec:HTPcomplete}

A computably enumerable set $C$ is said to be
\emph{$1$-complete} if every c.e.\ set $D$ is $1$-reducible to $C$,
written $D\leq_1 C$.  By definition this means that
for each $D$, there is a computable
total injective function $h:\omega\to\omega$ such that
$$ (\forall x\in\omega)~[x\in D \iff h(x)\in C].$$
The function $h$ is called a \emph{$1$-reduction}.  Of course,
the Halting Problem $\emptyset'$ is $1$-complete, and so
$1$-completeness of an arbitrary c.e.\ set $C$ is equivalent
to the statement $\emptyset'\leq_1 C$.  More generally,
for any set $A\subseteq\omega$ at all, the jump $A'$
is $1$-complete among sets computably enumerable in $A$,
in exactly the same sense.  Details appear in \cite[\S III.2]{S87}.

For a subset $W\subseteq\P$ of the set $\P$ of all prime numbers,
we define the subring $R_W=\Z[W^{-1}]$ in which the primes
with multiplicative inverses are precisely those in $W$.  The HTP operator maps $W$
to $HTP(R_W)$, as detailed in \cite{KM19}, and this set is clearly
computably enumerable relative to $W$, since a $W$-oracle allows one
to list out all rational numbers in $R_W$ and search for solutions to
polynomial equations.  Therefore we automatically have
$HTP(R_W)\leq_1 W'$.  In the case $W=\emptyset$, so that
$R_W=\Z$, the Matiyasevich-Davis-Putnam-Robinson result
shows that the reverse reduction also holds:  $\emptyset'\leq_1 HTP(R_{\emptyset})$.
This means that $\emptyset$ is \emph{HTP-complete}, according
to our new definition: its HTP is as complicated as possible.

\begin{defn}
\label{defn:HTPcomplete}
A set $W\subseteq\P$ of prime numbers is said to be \emph{HTP-complete}
if every set $V$ that is $W$-computably enumerable satisfies $V\leq_1 HTP(R_W)$.
Equivalently, $W$ is HTP-complete if and only if $W'\leq_1 HTP(R_W)$.
This is also equivalent to requiring $W'\equiv_1 HTP(R_W)$, or, by Myhill's Theorem,
to demanding that $W'$ and $HTP(R_W)$ be computably isomorphic.

It is also natural to say that $W$ is \emph{diophantine-complete} if every $V$ c.e.\ in $W$
is diophantine in the ring $R_W=\Z[W^{-1}]$.  However, the only sets $W$
currently known to have this property are the finite sets.
\end{defn}

In contrast to $\emptyset$, it is unknown whether $\P$ is HTP-complete,
since $R_{\P}=\Q$.  However, there is a broad result from \cite{KM19}
(see \cite[Cor.\ 3.3]{KM19} and the preceding remarks there), building on theorems
of Jockusch and Kurtz \cite{J81,K81}.  It will appear again in this article,
generalized as Theorem \ref{thm:JK}, with a proof.

\begin{thm}[see \cite{KM19}]
\label{thm:notHTPcomplete}
The set of all HTP-complete subsets of $\P$ is meager within the power
set of $\P$ and has Lebesgue measure $0$ there.
\end{thm}

This implies that the MRDP result for $\Z$ (that is, for $W=\emptyset$)
is an anomaly:  most subrings of $\Q$ do not have such
strong HTP's.  Of course, $\emptyset$ is hardly a representative element
of the power set of $\P$, so it is not surprising that it acts strangely.
Likewise, it would not be surprising if $\Q$ (that is, $R_{\P}$) did the same.
Nevertheless, it raises the question of just how many subsets of $\P$
are HTP-complete, and provides the initial answer: ``very few,''
in terms of both Lebesgue measure and Baire category.
In the next section, we balance this by showing that HTP-complete sets,
although meager, are ubiquitous:  they appear in every Turing degree,
and therefore there must be continuum-many of them.

\section{Building HTP-complete Sets}
\label{sec:construction}

\begin{thm}
\label{thm:construction}
Every Turing degree contains an HTP-complete set $V$.
Indeed, there is a uniform procedure which, given any set $C$,
computes an HTP-complete set $V\equiv_T C$.
\end{thm}
To be clear:  our procedure works uniformly on $C$, but
given two distinct Turing-equivalent sets $C$, it
will generally compute distinct sets $V$.  Indeed, one can make the
uniform procedure one-to-one, so that this is always the case.
It follows that there are countably many HTP-complete sets
Turing-equivalent to the given $C$.  Thus every Turing degree contains
infinitely many HTP-complete sets.

\begin{pf}
The following procedure,
using the oracle $C$, computes the required set $V$,
as we prove after giving the construction.  The construction
will begin with $V_0=\P$ and will delete elements
from $V$ at various stages, so that $V=\cap_s V_s$ is clearly
$\Pi^C_1$.  Afterwards we will argue that in fact $V\leq_T C$,
and then that $C\leq_T V$.

The requirement $\R_e$ demands that
$$ g_e \in HTP(R_V) \iff \Phi_e^C(e)\converges,$$
where $g_e$ is a polynomial we will define below.
For now it is acceptable to let the following polynomial
stand in for $g_e$, using the $e$-th odd prime $q_e$:
$$
f_e(X,Y,\ldots) =(X^2+q_eY^2-1)^2+\left(Y\left(1+\sum_{i=1}^4 Z_i^2
\right)-\left(1+\sum_{i=1}^4 W_i^2
\right)\right)^2
$$
whose solutions correspond to those pairs of rational numbers $(x,y)$
with $x^2+q_ey^2=1$ and $y>0$.  In fact, if $\Phi_e^C(e)\diverges$,
then this $f_e$ may have solutions in $\Q$, but will have no solutions
in a particular semilocal subring of $\Q$ determined in advance by the construction;
whereas if $\Phi_e^C(e)\converges$, then it will have solutions in every
semilocal subring.  This gives us a finite amount of wiggle room, enough
for the following finite-injury construction, and at the end we will replace
$f_e$ by a $g_e$ appropriate to the semilocal subring.

At each stage $s>\la e,0\ra$, $\R_e$ may protect various $q_e$-appropriate
primes $p_{e,t}$ from being removed from $V$.
If it ever sees $\Phi_{e}^C(e)$ converge, it will begin protecting
primes, and will protect them from then on, 
unless injured by a higher-priority requirement.  As long as $\Phi_{e,s}^C(e)$
diverges, its strategy is to remove from $V$ all $q_e$-appropriate
primes it can (but only finitely many at each stage).

We start by making a (uniformly computable) list of the primes
which each $\R_e$ will be allowed to protect.  Writing $p_{e,-1}=e$
for convenience, we define, for each $e\geq 0$ and $t\geq 0$:
$$  p_{e,t} = \min\set{\text{primes~}p>p_{e,t-1}}{
(\forall i \leq e+t)~[p\text{~is $q_i$-appropriate}\iff i=e]}.$$
Corollary \ref{cor:qpolys} shows that this set is nonempty,
and Lemma \ref{lemma:qpolys} shows that it is decidable uniformly in $e$ and $t$.
Thus $p_{e,0}$ is $q_i$-inappropriate for all $i<e$, so that
$\R_e$ will avoid any conflict with higher-priority requirements $\R_i$
which might need to remove $q_i$-appropriate primes from $V$.
The next prime
$p_{e,1}$ has all these properties and is also $q_{e+1}$-inappropriate,
so that if $\R_e$ comes to protect this prime, it will not injure
$\R_{e+1}$ by doing so.  As $\R_e$ protects increasingly larger primes,
it respects more and more requirements of lower priority.

At the stage $s+1$, we are given $V_s$, and we find the least
prime of the form $p_{e,t}$ that has not been considered at
any previous stage.  At this stage, we consider this prime,
fixing the $e$ and $t$ thus determined and writing $s'$
for the (earlier) stage at which $p_{e,t-1}$ was considered
(or $s'=0$ if $t=0$).  We compute $\Phi_{e,s}^C(e)$.
If this computation converges, then $\R_e$ continues to protect
all primes it protected at stage $s'$ (except any that may have
been removed from $V$ by higher-priority requirements in the interim),
and also protects $p_{e,t}$.  If it diverges, then $\R_e$ does not protect
any primes at this stage,
and deletes from $V_{s+1}$ those (finitely many) primes $p\in V_s$ satisfying:
\begin{itemize}
\item
$e < p < p_{e,t}$; and
\item
$p$ is $q_e$-appropriate; and
\item
$(\forall i<e)~\R_i$ is not protecting $p$ at this stage.
\end{itemize}
Thus $\R_e$ takes another step towards removing all $q_e$-appropriate primes from $V$,
since it still appears that $\Phi_e^C(e)$ diverges.
This completes the stage, and we set $V=\cap_s V_s$.

We remark that in this construction, if a prime $p$ is ever
protected by a requirement $\R_e$, then $p=p_{e,t}$ for some $t$
and no requirement will ever remove $p$ from $V$.
$\R_e$ will not: its computation must have converged
in order for it to have protected $p$ in the first place,
and so it will continue to protect $p$.  Moreover, its protection
stops any lower-priority $\R_j$ from removing $p$ from $V$.
Finally, higher-priority $\R_i$'s will not remove $p$ from $V$
because to have been chosen as $p_{e,t}$, $p$ must have
been $q_i$-inappropriate for each such $i$.

We now claim that $V$ is $C$-computable.  For a given prime $p$,
only requirements $\R_e$ with $e\leq p$ are ever allowed to
remove $p$ from $V$.  For each $e\leq p$, we can compute the
least number $t_e$ for which $p<p_{e,t_e}$.  If this $\R_e$ ever removes
$p$ from $V$, it must do so by the stage $s_e$ at which $p_{e,t_e}$ is considered;
the only reasons why it would not have removed $p$ from $V$
by this stage are either that $p$ is $q_e$-inappropriate or that
$\Phi^C_e(e)$ converged before stage $s_e$, or that a higher-priority
$\R_i$ is protecting $p$ at that stage.
In each of these cases, $\R_e$ will never remove $p$ from $V$.
So, by computing the maximum such stage $s=\max_{e\leq p}s_e$
and running the construction up to that stage $s$ (using a
$C$-oracle), we can check whether $p\in V$ or not.  Thus $V\leq_T C$.

Next we claim that every $\R_e$ is satisfied.  Here it is necessary
to define the specific polynomial $g_e$ to be used.  Let $f_e$
be as above, and let $g_e$ be derived from $f_e$ using Proposition
\ref{prop:semilocal}, so that, for every tuple $ (x,y,\zvec,\wvec)\in\Q^{10}$,
$$g_e(x,y,\zvec,\wvec)=0 \iff
f_e(x,y,\zvec,\wvec)=0~\&~(x,y,\zvec,\wvec)\in Q_e^{10},$$
where $Q_e$ is the semilocal subring of $\Q$ in which all primes
are inverted except those of the form
$p_{j,t}$ with $j+t \leq e$.  We claim that the map $e\mapsto g_e$
will be a $1$-reduction from $C'$ to $HTP(R_V)$.  Clearly this map is
injective, since each $f_e$ used a different coefficient $q_e$.
Moreover, it is computable, because the $C$-oracle is not invoked
in the definition of the primes $p_{e,t}$, nor for defining $g_e$.
So it remains to see that $e\in C'$ just if $g_e\in HTP(R_V)$.

Suppose first that $e\in C'$, and fix the least
stage $s+1$ at which we consider a prime of the form $p_{e,t}$
and for which $\Phi_{e,s}^C(e)\converges$.  From this stage on,
each prime $p_{e,t'}$ with $t'\geq t$ will be protected by $\R_e$,
starting at the stage at which it is considered.
Since it was chosen to be $q_i$-inappropriate for all $i<e$,
and since no lower-priority $\R_j$ can remove from $V$
a prime protected by $\R_e$, each such $p_{e,t'}$ lies in $V$
and gives a solution to $f_e$ in $R_V$.  As these primes
$p_{e,t'}$ are arbitrarily large, $g_e$ must lie in $HTP(R_V)$.

On the other hand, if $e\notin C'$, then $\R_e$ acts to remove
primes from $V$ at each stage $s+1$ at which any prime
$p_{e,t}$ is considered.  Therefore, it ultimately removes
from $V$ every $q_e$-appropriate prime $p>e$ except
those which are protected by higher-priority
requirements $\R_j$.  However, each prime $p_{j,t}$ protected
by $\R_j$ was chosen to be $q_i$-inappropriate for all $i\leq j+t$
except for $i=j$.  In particular, $p_{j,t}$ is $q_e$-inappropriate
whenever $e\leq j+t$; and if $e>j+t$, then $\frac1{p_{j,t}}\notin Q_e$.
Therefore, no $q_e$-appropriate prime is inverted
in the subring $(R_V\cap Q_e)$, and so $g_e\notin HTP(R_V)$.

This shows that $C'\leq_1 HTP(R_V)$, via the map $e\mapsto g_e$.
It follows that $C'\leq_1 V'$ and therefore $C\leq_T V$ by the Jump
Theorem \cite[Thm.\ III.2.3]{S87}.  On the other hand, with $V\leq_T C$,
we have $V'\leq_1 C'$, and of course $HTP(R_V)\leq_1 V'$, so $V$ is HTP-complete.
\qed\end{pf}

\begin{cor}
\label{cor:all}
For every Turing degree $\bfd\geq_T\bfz'$, there is a subring of $\Q$
for which Hilbert's Tenth Problem has Turing degree $\bfd$.
\end{cor}
\begin{pf}
This follows from Theorem \ref{thm:construction} along with the surjectivity
of the jump operator above $\bfz'$, which was first established by Friedberg
in \cite{F57}.
\qed\end{pf}

Many readers will recall \cite[Theorem 1.3]{P03} of Poonen,
described in detail in 
\cite[Ch.\ 12]{S07}.
It gives disjoint infinite decidable sets $T_1$ and $T_2$ of primes,
both of asymptotic density $0$, such that for every $W\supseteq T_1$
disjoint from $T_2$, $R_W$ admits a diophantine model of arithmetic
on the positive integers, making $\emptyset'\leq_1 HTP(R_W)$.
Unfortunately, the sets $W\subseteq\P$ containing all of $T_1$
(let alone disjoint from $T_2$) form a class that, while of
cardinality $2^\omega$, is meager and has measure $0$
in $2^{\P}$.  Thus this theorem cannot be combined with the results
in \cite{MCiE16} or \cite{MIMS}.  Its conclusion about
uncountability resembles Theorem \ref{thm:construction},
which showed that $\set{W}{W'\leq_1 HTP(R_W)}$ has size continuum,
despite being meager of measure $0$ .  If the goal is
to consider $\HTPQ$, then $\emptyset'\leq_1 HTP(R_W)$
seems just as relevant as $W'\leq_1 HTP(R_W)$,
and Poonen's stronger result about diophantine interpretation,
proven by an entirely different and deeper method than here,
could be the key to a final answer.  (Cf.\ Theorems 41 and 53 of \cite{MCiE16}.)

We also remark briefly that using the foregoing method with $C=\emptyset$,
one can also prove that there is a decidable subring $R_W\subseteq\Q$
for which $HTP(R_W)\equiv_1\emptyset'$.  Of course, the
Matiyasevich-Davis-Putnam-Robinson theorem already proved this result
for the far more challenging specific case $R=\Z$.
Still, the results here again suggest how a
computability-theoretic approach, using techniques such as
finite-injury constructions along with basic number theory,
can sometimes yield new and different proofs.
One continues to hope that those techniques, combined with a deeper use
of number theory than in this article, might accomplish more
than either discipline can achieve on its own.

Our method here does not appear to provide answers to any of the
questions raised in \cite[Remarks 3.20 \& 4.8]{EMPS15}
by Eisentr\"ager, Miller, Park, and Shlapentokh.
Those questions generally want the degree of $HTP(R_W)$
to be held down, so that $W'\not\leq_T HTP(R_W)$,
whereas the method of this section is appropriate for coding information
into $\HTP{R_W}$ and thus making its Turing degree large.

\section{Enumeration Operators}
\label{sec:enumops}

Theorem \ref{thm:notHTPcomplete} was proven with no use of the HTP-operator
specifically.  It used only the fact that $HTP$ is an \emph{enumeration operator}.
Recall that an \emph{enumeration} of a subset $A$ of $\omega$ is a subset $B$
of $\omega$ that, when viewed as a subset of $\omega^2$, projects onto $A$
via the projection $\pi_1$:
$$ A = \pi_1(B) = \set{x\in\omega}{(\exists y)~\la x,y\ra\in B}.$$
This is equivalent to various other definitions (often using functions).

\begin{defn}
\label{defn:enumop}
Let $G:2^\omega\to 2^\omega$ be a function.  $G$ is said to be an
\emph{enumeration operator} if there exists a Turing functional $\Gamma$
such that, for every $A\in 2^\omega$ and every enumeration $C$ of $A$,
$\Gamma^C$ is a total function from $\omega$ into $\{ 0,1\}$ and is the characteristic
function of an enumeration of $G(A)$.  (It is also natural to refer to $\Gamma$ itself
as the enumeration operator, but it confuses matters.  We will say here that $\Gamma$
\emph{represents} $G$.)

We write $B\leq_e A$, and say that $B$ is \emph{enumeration-reducible to $A$},
or \emph{$e$-reducible to $A$}, if there is an enumeration operator $G$
with $G(A)=B$.  This is equivalent to the usual definition, e.g.\ in \cite[\S 9.7]{R67}.
\end{defn}
It is immediate from the definition that if $\Gamma$ represents an enumeration operator
and $\pi_1(C_0)=\pi_1(C_1)$, then $\pi_1(\Gamma^{C_0}) = \pi_1(\Gamma^{C_1})$.
We note that other definitions of $e$-reducibility are standard in the literature, and
are readily shown to be equivalent to this one.  The essence is that there exists
a uniform procedure that accepts any enumeration of $A$ and uses it to
compute an enumeration of $G(A)$.

The jump operator $J$, mapping each $A$ to $A'$, is the prototype of the functions called
\emph{pseudojump operators} by Jocksuch and Shore in \cite{JS83,JS84},
whose output can be enumerated uniformly when we are given $A$ itself
(not just an enumeration) as the oracle.
\comment{
\begin{defn}
\label{defn:PJop}
A function $G:2^\omega\to 2^\omega$ is a \emph{pseudojump operator}
if there exists an index $e$ such that, for all $A\in 2^\omega$,
$$ G(A) = A\oplus W_e^A = A\oplus \dom{\Phi_e^A},$$
where $W_e^A=\dom{\Phi_e^A}$,
using the standard list of Turing functionals $\Phi_e$.
\end{defn}
}
However, the jump operator is not an enumeration operator.
To see this, notice that if it were, then $\emptyset''=J(\emptyset')$ would also be 
computably enumerable, since we could run the representation $\Gamma$
on a computable enumeration of $\emptyset'$ to get a computable enumeration
of $J(\emptyset')$.  For a better understanding of the failure of the jump to be
an enumeration operator, consider a functional $\Phi_e$ for which, for all $x$,
$$ \Phi_e^A(x) =\left\{\begin{array}{cl} 0, & \text{if~}17\notin A;\\
\diverges, & \text{if~}17\in A.\end{array}\right.$$
Now if some functional $\Gamma$ represented the jump (as an enumeration operator),
then with $A=\emptyset$ we would have $\Gamma^{\emptyset}(e)=1$, as $\emptyset$ itself
is an enumeration of $\emptyset$.  But if $u$ is the use of this computation, then one readily
can create an enumeration $C$ of an arbitrary $B\subseteq\omega$ with $C\res u=0^u$,
and $\Gamma^C(e)$ would have to equal $1$ for each such $C$, by the Use Principle.
Hence $\Gamma$ either fails to be an enumeration operator, or else fails to compute
the jump, because many sets $A$ (indeed a class of measure $\frac12$) have $e\notin A'$.

The next result generalizes Theorem \ref{thm:notHTPcomplete}, and we now give
a proof, by exactly the same means as in \cite[Cor.\ 3.3]{KM19}.
\begin{thm}
\label{thm:JK}
For every enumeration operator $E$, the collection $\set{A\in 2^\omega}{A'\leq_1 E(A)}$
is meager and has measure $0$.
\end{thm}
\begin{pf}
With $E$ fixed, we show that $A'\not\leq_1 E(A)$ for every set $A$ such that,
for some set $B<_T A$, $A$ is $B$-computably enumerable.  Indeed,
$E(A)$ must then also be $B$-c.e., so $E(A)\leq_1 B'$.  However, with
$A\not\leq_T B$, we have $A'\not\leq_1 B'$, by the Jump Theorem
(see e.g.\ \cite[Thm.\ III.2.3]{S87}).  It would now contradict
the transitivity of $1$-reducibility to have $A'\leq_1 E(A)$.

By results of Jockusch and his student (at the time) Kurtz in \cite{J81,K81},
the class of \emph{relatively c.e.\ sets}, i.e., those $A$ for which a $B$ exists as
described above, is a comeager class of measure $1$.  The theorem follows.
\qed\end{pf}

On the other hand, it is quite possible for $\set{A\in 2^\omega}{A'\leq_T E(A)}$
to be comeager and to have measure $1$.  Indeed, the enumeration operator
mapping $A$ to $(\emptyset'\oplus A)$ has this property:  it is well-known
that the class $GL_1$ of \emph{generalized-low} sets, i.e., those satisfying
$A'\leq_T \emptyset'\oplus A$, is comeager and has full measure.  Here we
focus on the possibility of computing $A'$ uniformly (via a single Turing functional)
from $E(A)$.

\begin{thm}
\label{thm:epsilon}
For every Turing functional $\Psi$ and every enumeration operator $E$,
$\mu(\set{A\in 2^\omega}{\chi_{A'}=\Psi^{E(A)}})<1$.
\end{thm}

\begin{cor}
\label{cor:epsilon}
For every Turing operator $\Phi$, there exists a set $\mathcal S$ of positive measure
such that, for all $W\in\mathcal S$,
$$ \Phi^{HTP(R_W)} \neq \chi_{W'}.$$
\end{cor}

\begin{pftitle}{Proof of Theorem \ref{thm:epsilon}}
We fix an index $e$ for a Turing functional defined as follows:
$$ \Phi_e^A(x) = \left\{\begin{array}{cl}
0, & \text{if~}(\exists m>1)~\{ m+1,m+2,\ldots,2m\}\cap A=\emptyset;\\
\diverges, & \text{otherwise}.
\end{array}\right.$$
Thus the measure of the set of those $A$ with $e\in A'$ is at most $\frac12$
(in fact, somewhat less than $\frac12$ because of overlaps) and certainly positive.
Suppose that, on a set of measure $1$, $\Psi^{E(A)}=A'$ (that is, $\Psi^{E(A)}$ computes the
characteristic function of $A'$).  Then $\mu(\set{A}{\Psi^{E(A)}(e)\converges=0})>0$.
By the countable additivity of Lebesgue measure, there must exist a specific $\sigma\in 2^{<\omega}$
such that $\Psi^{\sigma}(e)\converges=0$ and such that
$\mu(\set{A}{\sigma\initial E(A)})>0$.  Indeed, since the relation $\sigma^{-1}(1)\subseteq E(A)$
is always created by a finite subset of $A$, there must then exist a finite
set $S_0$ such that $\sigma\initial E(S_0)$ and
$\mu(\set{A}{S_0\subseteq A~\&~\sigma\initial E(A)})>0$, since there are only
countably many finite subsets $S$ of $\omega$.
(To avoid confusion, in this proof we write $\sigma\initial E(A)$ to mean that $\sigma$
is an initial segment of $E(A)\in 2^\omega$, and $S\subseteq A$ to mean
simply that $S$ is a subset of $A$, not necessarily an initial segment.)

We fix such a $\sigma$ and such an $S_0$, and choose an integer
$m >\max(S_0)$ (with $m>1$ as well).  Let $\W = \set{A}{S_0\subseteq A~\&~\sigma\initial E(A)}$,
which is thus guaranteed to have positive measure.  Now consider the class
$$ \V = \set{B\in 2^\omega}{(\exists A\in\W)~B=A-\{ m+1,m+2,\ldots,2m\} }.$$
For every such $B$, $m$ witnesses that $e\in B'$, according to our definition
of $\Phi_e$.  (In contrast, only measure-$0$-many $A\in\W$ lie in $\V$, since $e\notin A'$.)
Moreover, since $m>\max(S_0)$, all $B\in\V$ have $S_0\subseteq B$
and thus have $E(S_0)\subseteq E(B)$, since enumeration operators are clearly
monotone under $\subseteq$.  On the other hand, each $B\in\V$ has a corresponding $A\in\W$
for which $B\subseteq A$, so that $E(B)\subseteq E(A)$.  Together these yield
$\sigma\initial E(B)$, since $E(S_0)$ and $E(A)$ agree up to $|\sigma|$.
But now, for every $B\in\V$, we have $\Psi^{E(B)}(e)\converges=\Psi^{\sigma}(e)=0$,
even though $e\in B'$.

It remains to show that $\V$ has positive measure.  Suppose $\set{\U_{\tau_i}}{i\in\omega}$
is a cover of $\V$ by basic open subsets $\U_{\tau_i}=\set{C}{\tau_i\initial C}$ of
Cantor space, and suppose that this cover has total Lebesgue measure $r$.
By the definition of $\V$, we may assume that $\tau_i(n)=0$ for all $n\in\{ m+1,\ldots,2m\}$
and all $i$ with $n<|\tau_i|$.  Now, for each of the $(2^m-1)$-many binary strings
$\rho$ of length $m$ (excluding the zero string $0^m$), let $\T_{i,\rho}=\U_{\tau_{i,\rho}}$,
where
$$ \tau_{i,\rho}(n) = \left\{\begin{array}{cl}
\rho(n-(m+1)), & \text{if~}m+1\leq n \leq 2m~\&~n<|\tau_i|;\\
\tau_i(n), & \text{otherwise.}
\end{array}\right.$$
That is, $\tau_{i,\rho}$ is the same as $\tau_i$, except that the portion from $(m+1)$
up to $2m$, which was all zeroes in $\tau_i$, is replaced by the (nonzero) string $\rho$.
Thus $\mu(\U_{\tau_i})=\mu(\T_{{i,\rho}})$ for all $i$ and $\rho$.  Since the sets $\U_{\tau_i}$
form a cover of $\V$, the definition of $\V$ shows that the sets $\T_{i,\rho}$ form
an cover of $\W$ by basic open sets in Cantor space, so that their total measure
is $\geq \mu(\W) >0$.  Also, for any two
distinct $\rho$ (and the same $i$), the strings $\tau_{i,\rho}$ are distinct;
whereas for distinct $i$ and the same $\rho$, the overlap between strings $\tau_{i,\rho}$
is equal in measure to the overlap between the corresponding $\tau_i$.
It follows that
$$ \mu(\W) \leq \mu\left(\bigcup_{i\in\omega}\bigcup_{\text{~nonzero~}\rho\in 2^m} \T_{i,\rho}\right)
= (2^m-1)\cdot\mu\left(\bigcup_i \U_{\tau_i}\right).$$
Therefore, this open cover $\{\U_{\tau_i}\}$ of $\V$ has Lebesgue measure
at least $\frac{\mu(\W)}{2^m-1}$, and this positive lower bound is independent
of the choice of cover of $\V$.  So $\mu(\V)$ is positive as well, and we saw above
that $\Psi^{E(B)}(e)\converges\neq B'(e)$ for all $B\in\V$.
\qed\end{pftitle}

It remains possible, therefore, that $W'\leq_T HTP(R_W)$ might hold for
measure-$1$-many sets $W$, but if so, it requires infinitely many Turing functionals
to establish this fact.  Similarly, the reduction $A'\leq_T \emptyset'\oplus A$
can be established on a set of measure $(1-\epsilon)$ by a single functional $\Phi$
(for arbitrarily small $\ep>0$), but countably many functionals are required
to show that it holds on a set of measure $1$.  (In that case, the countably many
functionals can be produced uniformly in the rational number $\ep>0$.)


\section{Existential Definability of $\Z$}
\label{sec:definability}

Existential definability of a subset $S$ of $\Q$ (in the usual
model-theoretic notion, i.e., defining a unary relation
on the field $\Q$ whose elements are precisely the elements of $S$)
is equivalent to $S$ being \emph{diophantine} in the ring $\Q$.
This means that, for some $n$, $S$ is defined by a single polynomial
$f\in\Q[X,Y_1,\ldots,Y_n]$ as follows:
$$ (\forall r\in\Q)~[r\in S \iff (\exists\yvec\in \Q^n)~f(r,\yvec)=0].$$
All more complicated existential definitions can be boiled down
to definitions of this form.

It is unknown whether the set $\Z$ is existentially definable
in the field $\Q$.  Julia Robinson gave the first definition
of $\Z$ in $\Q$, in \cite{R49}.  That definition was $\Pi_4$.
Significant subequent work has reduced the complexity
of such definitions:  Poonen \cite{P09} gave a $\Pi_2$ definition,
and then Koenigsmann \cite{K15} gave a $\Pi_1$ (that is,
purely universal) definition.  Thus we seem to be getting
closer to an existential definition.  However, there are number-theoretic
conjectures, notably by Mazur, that would imply the existential undefinability
of $\Z$ in $\Q$.

An existential definition of $\Z$ in $\Q$ would imply
$HTP(\Z)\leq_1 HTP(\Q)$, and hence $\emptyset'\leq_1 HTP(\Q)$,
so it is highly relevant to this article.  However, our purpose in this section is
to investigate other possible consequences of $\exists$-definability
of $\Z$ in $\Q$.  The main point is that, if any of these consequences
should be shown not to hold, it would follow that $\Z$ is not
diophantine in $\Q$.

From an existential definition of $\Z$ within the field $\Q$,
we would immediately get a stronger result.
\begin{lemma}
\label{lemma:allsubrings}
If $\Z$ has an existential definition in $\Q$, then indeed there is a polynomial
$h\in\Z[X,Y_1,\ldots,Y_k]$ such that, for all $x\in\Q$,
$$ x\in\Z \iff (\exists \yvec\in\Q^k)~h(x,\yvec)=0
\iff (\exists \yvec\in\Z^k)~h(x,\yvec)=0.$$
Thus the formula $(\exists Y_1\cdots\exists Y_k)h(X,\Yvec)=0$
would define $\Z$ not only in $\Q$, but also in every subring of $\Q$.
Likewise, every c.e.\ set would have an existential definition
that holds in every subring of $\Q$.
\end{lemma}
\begin{pf}
Assume that the formula $(\exists Z_1,\ldots,Z_j)~g(X,\Zvec)=0$
defines $\Z$ in $\Q$, with $g$ of total degree $d$.  Define $h(X,\Yvec,\Tvec)$
to be the polynomial
$$ g\left(X,\frac{Y_1}{1+T_1^2+\cdots+T_4^2},\ldots,\frac{Y_j}{1+T_1^2+\cdots+T_4^2}\right)
\cdot (1+T_1^2+\cdots+T_4^2)^{d}$$
Now if $x\in\Z$, then there is $\zvec\in\Q^j$ with $g(x,\zvec)=0$.
Taking a positive common denominator $v\in\Z_{>0}$ of the rationals $z_i$,
use the Four Squares Theorem to write $v-1=t_1^2+t_2^2+t_3^2+t_4^2$
with all $t_i\in\Z$ and let $y_i=vz_i$.
Then $\left(\yvec,\tvec~\!\right)$ is a solution to $h(x,\Yvec,\Tvec)=0$
in $\Z$, hence in every subring of $\Q$.

Conversely, for any $\left(x,\yvec,\tvec~\!\right)$ in a subring of $\Q$ with $h\left(x,\yvec,\tvec~\!\right)=0$,
setting $z_i=\frac{y_i}{1+t_1^2+t_2^2+t_3^2+t_4^2}$ gives
$g(x,\zvec)=0$ with all $z_i\in\Q$, so $x\in\Z$.
\comment{
This is simply a matter of writing the rational variables as quotients
of integer variables, clearing denominators, and applying the
Four Squares Theorem.  Once $\Z$ is defined thus, it is easy to
translate any existential definition (of an arbitrary c.e.\ set) in $\Z$
into an existential formula that either holds in $\Z$ or fails in all of $\Q$.
}
\qed\end{pf}

\subsection{Preservation of $m$-reductions}
\label{subsec:leqm}

It was seen in \cite{KM19} that the HTP operator can fail
to preserve Turing reductions, and indeed that it can sometimes
reverse them:  it is possible to have $V <_T W$, yet
$HTP(R_W) <_T HTP(R_V)$, with strictness in both
relations.  (This result is \cite[Corollary 5.3]{KM19}.)  Whether the same operator
must respect the stronger notion of $m$-reducibility remains
an open question.  Here we connect that question to the existential
definability of $\Z$ in $\Q$, first giving the relevant definitions.
\begin{defn}
\label{defn:leqm}
For subsets $A,B\subseteq\omega$, a computable total function $F:\omega\to\omega$
is an \emph{$m$-reduction} from $A$ to $B$ if it satisfies
$$ (\forall x\in\omega)~[x\in A \iff F(x)\in B].$$
A \emph{$1$-reduction} is just an $m$-reduction which is also
one-to-one (as opposed to \emph{many-to-one}, whence the terminology).
We write $A\leq_1 B$ and $A\leq_m B$ to denote the existence
of $1$-reductions and $m$-reductions, respectively.  Clearly these
are both partial preorders on the power set of $\omega$.
\end{defn}

The reader may wonder why the distinction is made between
$m$- and $1$-reducibility.
There do exist sets $A$ and $B$ with $A\leq_m B$ but $A\not\leq_1 B$,
and they can be chosen to be infinite and coinfinite (thus avoiding
the simple situation where $1\leq |B| < |A| <\infty$).  Nevertheless, in computability
theory, $1$-reducibility is regarded as nearly equivalent to $m$-reducibility.
Our first lemma suggests that this seems to hold here as well.

\begin{lemma}
\label{lemma:1vsm}
For sets $A\subseteq\omega$ and $W\subseteq\P$, we have $A \leq_m HTP(R_W)$
if and only if $A \leq_1 HTP(R_W)$.
\end{lemma}
\begin{pf}
For the nontrivial direction, let $G$ be an $m$-reduction.  Then
each value $G(n)$ is a polynomial in $\Z[X_1,X_2,\ldots]$, say,
and we simply define:
$$ F(n) = (G(n))^2 + (X_0)^{2n}.$$
The polynomial $F(n)$ (from $\Z[X_0,X_1,\ldots]$) has a solution in $R_W$ just if $G(n)$ does,
and the exponent $2n$ makes $F$ injective.
\qed\end{pf}
\begin{cor}
\label{cor:1vsm}
If the HTP operator respects $m$-reductions, then it respects $1$-reductions.
\qed\end{cor}
Nevertheless, there is an important reason to distinguish between
$1$- and $m$-reducibility, as seen in the following theorem.
\begin{thm}
\label{thm:equivalence}
Each of the following implies the next.
\begin{enumerate}
\item
$Z$ is existentially definable in the field $\Q$.
\item
The HTP operator respects $m$-reducibility
(i.e., if $V\leq_m W$, then $ HTP(R_V)\leq_m HTP(R_W)$).
\item
$\emptyset'\leq_1\HTPQ$.
\end{enumerate}
\end{thm}
In contrast, we do not know whether (3) follows from the assumption
that HTP preserves $1$-reductions.
\begin{pf}
We first show that (2) implies (3).  Consider $V=\{ 3\}$ and $W=\P-\{3\}$.
Clearly $V\leq_m W$:  just let $F(3)=5$ and $F(p)=3$ for all $p\neq 3$.
(This would work for any nonempty finite $V$ and any proper cofinite $W$, of course.)
By (2), we get $HTP(\Z[\frac13])\leq_m HTP(R_{\P-\{3\}})$.
But Julia Robinson showed that $\emptyset'\leq_1 HTP(\Z[\frac13])$
(and likewise for all finitely generated subrings of $\Q$), whereas
$HTP(R_{\P-\{3\}})\leq_1\HTPQ$ by Corollary \ref{cor:semilocal},
proving (3).

Next we assume (1) and prove (2).
With an $m$-reduction from $V$ to $W$, we can readily compute
an $m$-reduction from $R_V$ to $R_W$:  that is, a computable, total,
function $G$ with
$$ (\forall q\in\Q)~[ q\in R_V \iff G(q)\in R_W].$$
$\exists$-definability of $\Z$ implies that every c.e.\ set, and in particular the graph of $G$,
is diophantine in $\Q$, so by Lemma \ref{lemma:allsubrings}
we have a polynomial $g$ such that,
for all $q,r\in\Q$:
\begin{align*}
G(q)=r &\iff g(q,r,Z_1,\ldots,Z_m)\in\HTPZ\\
&\iff  g(q,r,Z_1,\ldots,Z_m)\in\HTPQ.
\end{align*}
Thus the following holds of every $f\in\Z[X_0,\ldots,X_{k-1}]$:
\begin{align*}
f\in HTP(R_V) &\iff (\exists\qvec\in (R_V)^k)~f(\qvec)=0\\
\iff& (\exists\qvec\in\Q^k)(\exists\rvec\in (R_W)^k)~[f(\qvec)=0~\&~(\forall i<k)~G(q_i)=r_i]\\
\iff& (\exists\qvec\in\Q^k)(\exists\rvec\in (R_W)^k)\\
&~[f(\qvec)=0~\&~(\forall i<k)~g(q_i,r_i,Z_{i1},\ldots,Z_{im})\in HTP(R_W)]\\
\iff& (\exists \svec,d,\rvec,z_{01},\ldots,z_{km} \in R_W)~\left[f\left(\frac{s_1}d,\ldots,\frac{s_k}d\right)=0~\&\right.\\
&~\left.\&~(\forall i<k)~g\left(\frac{s_i}d,r_i,z_{i1},\ldots,z_{im}\right)=0~\&~d\neq 0\right]
\end{align*}
\comment{
\begin{align*}
f\in HTP(R_V) &\iff (\exists\qvec\in (R_V)^k)~f(\qvec)=0\\
\iff& (\exists\qvec\in\Q^k)(\exists\rvec\in (R_W)^k)~[f(\qvec)=0~\&~(\forall i<k)~G(q_i)=r_i]\\
\iff& (\exists\qvec\in\Q^k)(\exists\rvec\in (R_W)^k)\\
&~[f(\qvec)=0~\&~(\forall i<k)~g(q_i,r_i,Z_{i1},\ldots,Z_{im})\in\HTPQ]\\
\iff& (\exists \svec,d,\rvec,\tvec_1,\ldots,\tvec_m,u \in R_W)~[f\left(\frac{s_1}d,\ldots,\frac{s_k}d\right)=0~\&\\
&~\&~(\forall i<k)~g\left(\frac{s_i}d,r_i,\frac{t_{i1}}u,\ldots,\frac{t_{i_m}}u\right)=0~\&~du\neq 0]\\
\iff& \text{the polynomial in the line above lies in~}HTP(R_W)
\end{align*}
}
Since the equations (and the inequation) in the second-to-last line
can all be collected into a single polynomial equation with the $d$'s cleared
from the denominators, we have computed (from $f$) a single polynomial
which lies in $HTP(R_W)$ just if $f$ itself lies in $HTP(R_V)$.
\qed\end{pf}

\subsection{Boundary Sets of Polynomials}
\label{subsec:boundary}

The key to our use of the polynomials $f_e$ built using $(X^2+q_eY^2-1)$ in Theorem
\ref{thm:construction}, and also in the results in \cite{KM19},
was that, once we built $f_e$ and thus ruled out the trivial solutions, they have
nonempty \emph{boundary sets}, according to the following definition.

\begin{defn}
\label{defn:ABC}
For a polynomial $f\in\Z[X_1,X_2,\ldots]$, write:
\begin{itemize}
\item
$\A(F)=\set{W\in 2^{\P}}{f\in HTP(R_W)}$;

\item
$\C(F)=\set{W\in 2^{\P}}{(\exists\text{~finite~}S_0\subseteq \overline{W})~f\notin HTP(R_{\P-S_0})}$;

\item
$\B(F)=2^{\P}-\A(f)-\C(f)$; the \emph{boundary set} of $f$.
\end{itemize}
With $\mu$ as the Lebesgue measure on $2^{\P}$, we also write
$\alpha(f)=\mu(\A(f))$, $\beta(f)=\mu(\B(f))$, and $\gamma(f)=\mu(\C(f))$.
\end{defn}
The Cantor space $2^{\P}$ is equipped with the usual topology.
Here $\A(f)$ is always an open set, since each solution to $f$ requires
only that a certain finite set of primes be inverted in $R_W$.
$\C(f)$ is the interior of the complement of $\A(f)$, the set of subrings
where the non-invertibility of some finite set of primes rules out
the possibility of a solution to $f$.  Therefore, $\B(f)$ is indeed the
topological boundary of $\A(f)$, and contains those $W$ such that
$f$ has no solution in $R_W$, but such that, for every $n$, it is possible
to extend $W\res n$ to some set $V$ with $f\in HTP(R_V)$.
(In the phrase of Alexandra Shlapentokh, $f$ ``never loses hope''
of having a solution in $R_W$.)  Often $\B(f)$ is empty, but the
polynomials $f_e$ have nonempty boundary sets:  indeed
$\B(f_e)$ contains every subset of the set of $q_e$-inappropriate
primes, which Lemma \ref{lemma:qpolys} showed to be an infinite set.
This is what allowed our coding to work, in Theorem \ref{thm:construction}:
no matter how many primes we removed from $V$, there was always
some prime not yet removed which, if it stayed in $V$,
would cause $f_e$ to lie in $HTP(R_W)$.  So, no matter how long
$\Phi_e^C(e)$ might take to converge, we could always code
its convergence into $HTP(R_W)$ when and if we saw the computation halt.

On the other hand, the definitions of $\alpha$, $\beta$, and $\gamma$
suggested that we care about the measures of these sets, and here the
$f_e$ polynomials are not so impressive.  Indeed, $\alpha(f_e)$ is always $1$,
for every $e$, because the set of $q_e$-appropriate primes is infinite
and the inversion of any single element of that set will yield a solution
to $f_e$.  It remains an open question whether any polynomial $f$ at all
can have $\beta(f)>0$.  In this section we discuss the possible consequences
of an answer to this question.

The overall boundary set $\B$ is defined by:
$$ \B = \bigcup_{f\in\Z[X_1,X_2,\ldots]}\B(f).$$
Each $\B(f)$ is nowhere dense in $2^{\P}$, in the sense of Baire category,
and therefore $\B$ itself is meager.  This shows that there must exist
subrings of $\Q$ which lie in no boundary set $\B(f)$.  These are called
\emph{HTP-generic} subrings, and are studied in \cite{MCiE16, MIMS}.
As noted above, although the complement $\overline{\B}$ is comeager
and thus large in the sense of Baire category, it is unknown whether
its measure is $0$ or $1$, and even values between $0$ and $1$ have not
been ruled out.  We remark that, for an individual polynomial $f$,
we always have $\beta(f)<1$, because the only way to have $\alpha(f)=0$
is for the open set $\A(f)$ to be empty, in which case $\C(f)=2^{\P}$
and $\B(f)=\emptyset$.

\subsection{Noncomputable $\beta(f)$}
\label{subsec:Sigma1}

The next theorem will be superseded by Theorem \ref{thm:rprime},
but its proof is useful as an introduction to the proof of the latter theorem,
and so we present it in full here.

\begin{thm}
\label{thm:measure0}
If the boundary set $\B$ has measure $<1$, then there is no
existential definition of $\Z$ within the field $\Q$.
\end{thm}
\begin{pf}
We prove the contrapositive, by assuming that $\Z$ does have an
$\exists$-definition in $\Q$ and showing, for an arbitrary positive real
number $r < 1$ which is approximable from below, that there exists
a polynomial $f\in\Z[\Xvec]$ with $\alpha(f)=r$ and $\gamma(f)=0$.
(``Approximable from below'' means that the left Dedekind cut of $r$ is c.e.)
This will establish that the measure $\beta(f)=1-r$, proving the theorem,
since $\B(f)\subseteq\B$.

So fix such a number $r$, and let $q_0,q_1,\ldots$ be a computable,
strictly increasing sequence of positive rational numbers with $\lim_s q_s=r$.
Let $n_0$ be the least integer with $2^{-n_0}\leq q_0$, which is to say,
$1-2^{-n_0}\geq 1-q_0$.  Now define by recursion
$$ n_{k+1}=\min\set{n\in\N}{(1-2^{-n})\cdot (1-2^{-n_k})\cdots(1-2^{-n_0})\geq 1-q_{k+1}}.$$
With $q_{k+1}>q_k$, such an $n_{k+1}$ always exists (and must be positive,
since $q_{k+1} < 1$), and the sequence $\la n_i\ra_{i\in\N}$ is computable.
Moreover, by the minimality of each $n_k$,
$\prod_{k\geq 0} (1-2^{-n_k})=1-\lim_k q_k = 1-r$.

Next, let $x_0=p_0\cdot p_1\cdots,p_{n_0-1}$ be
the product of the first $n_0$ prime numbers. Then set
$x_{k+1}=p_{n_0+\cdots+n_k}\cdots p_{n_0+\cdots+n_{k+1}-1}$
to be the product of the next $n_{k+1}$ primes, for each $k$ in turn.
The set $D=\set{x_k}{k\in\N}$ is computably enumerable (indeed computable),
hence diophantine.  Since we are assuming that $\Z$ is $\exists$-definable
in $\Q$, there exists a polynomial $g\in\Z[X,Y_1,\ldots,Y_m]$ such that
\begin{align*}
D &= \set{x\in\Z}{g(x,Y_1,\ldots,Y_m)\in\HTPZ}\\
&=\set{x\in\Q}{g(x,Y_1,\ldots,Y_m)\in\HTPQ}.
\end{align*}
(This simply requires that we start with a polynomial which defines the set
$D=\set{x_k}{k\in\N}$ within $\Z$, and then apply Lemma
\ref{lemma:allsubrings} to transfer the definition of $\set{x_k}{k\in\N}$ to $\Q$.)
The $f(X,\Yvec,T)$ we desire is simply the sum
$$ (g(X,Y_1,\ldots,Y_m))^2
+(XT-1)^2.$$
We claim that this $f$ satisfies $\alpha(f)=r$ and $\gamma(f)=0$.

Notice first that every solution $(x,\yvec,t)$ to $f$ in $\Q$
must have $g=0$, hence has $x\in\Z$ and all $y_i\in\Z$.
But then $x=x_k$ for some $k$, by our choice of $g$, and $t=\frac1{x_k}$.
In order for this solution to lie in a subring $R$ of $\Q$, therefore,
that subring $R$ must contain multiplicative inverses of all the prime factors
$p_{n_0+\cdots+n_{k-1}},\ldots,p_{n_0+\cdots+n_k-1}$ of this $x_k$.
(Notice that this list contains exactly $n_k$ primes.)

Conversely, suppose that a subring $R$ does contain all these
primes (for some $k$).  Then it contains $t=\frac1{x_k}$, and
since $x_k\in\N$, there exist integers $y_1,\ldots,y_m$ which,
along with $x_k$ and $t$, form a solution to $f$ in $R$.

Therefore, the subrings in which $f$ has a solution are exactly
those in which, for some $k$, all of the $n_k$ prime factors
of $x_k$ have inverses.  For a single $k$, the measure of the set of such
subrings is $2^{-n_k}$.  Since all distinct $x_k$ have completely
distinct prime factors, the set of subrings containing no solution to $f$
therefore has measure
$$\prod_k (1-2^{-n_k}) = 1-r,$$
and so the set $\A(f)$ of subrings with solutions to $f$
has measure precisely equal to $r$.  That is, $\alpha(f)=r$.

Finally, it is clear that every semilocal subring $R$ of $\Q$ contains a solution of $f$.
Indeed, for some $k$, $R$ must contain inverses of all primes $\geq p_{n_0+\cdots+n_k}$,
so our analysis above yields a solution in $R$.  It follows that $\C(f)=\emptyset$,
so $\beta(f)=1-\alpha(f)-\gamma(f)=1-r$.
\comment{
Some additional work is required to complete the proof under the assumption
that the ring $\Z$ has a diophantine interpretation in the field $\Q$,
but the essential ideas are the same.
}
\qed\end{pf}

The proof of Theorem \ref{thm:measure0} actually showed more.
Assuming an $\exists$-definition of $\Z$ in $\Q$, we constructed
a polynomial $f$ with $\alpha(f)=r$ and $\beta(f)=1-r$, under the condition
that $r\in (0,1)$ be approximable from below.  In particular, this shows
that both $\alpha(f)$ and $\beta(f)$ can be noncomputable, since a real number
$r$ can be approximable from below without being approximable from above.
\begin{cor}
\label{cor:from above}
If $\Z$ has an existential definition in $\Q$, then for every real
number $r\in (0,1)$ which is approximable from below, then there
is a polynomial $f\in\Z[\Xvec]$ with $\alpha(f)=r$ and $\beta(f)=1-r$.
\qed\end{cor}

Finally, we remark that in the proof of Theorem \ref{thm:measure0},
it is possible to put an upper bound on the degrees of the polynomials
$f$ produced.  First of all, the polynomials $h$ and $j$ (and $(XT-1)$)
are all fixed independently of $r$, and hence so is the total degree $d$
of $h$.  Only $g$ depends on $r$:  $g$ was chosen to define (in $\Z$)
the set $D$ of products $x_k$ of primes, and the number of prime factors
of each $x_k$ depends on $r$.  However, by fixing a single polynomial
$G\in\Z[E,X,Y_1,\ldots,Y_k]$ which defines the Halting Problem in $\Z$,
we may then take our $g$ (for a given $r$) to be of the form $G(e,X,\Yvec)$
for some natural number $e$.  (In fact, the choice of $e$ can be made
effectively, once we know an index for the computable sequence
$\la q_k\ra_{k\in\N}$ of rationals approaching $r$ from below.)
Therefore, regardless of the value of $r$, the total degree of $g$
need never be more than that of the fixed polynomial $G$, and this
in turn puts a bound on the total degree of the $f$ we eventually produced.

\subsection{Reals of Greater Complexity}
\label{subsec:Sigma2}

Having seen in Subsection \ref{subsec:Sigma1} how to use
arbitrary c.e.\ sets, along with the assumption
of $\exists$-definability of $\Z$ in $\Q$, to build polynomials
$f$ with $\beta(f)$ noncomputable, we now enhance our
construction of the c.e.\ set, so as to make $\beta(f)$ have
even higher complexity.  From its definition, $\beta(f)=1-\alpha(f)-\gamma(f)$,
and $\alpha(f)$ must always be approximable from below,
while $\beta(f)$ must be $\HTPQ$-approximable from below,
hence $\HP$-approximable from below.
We will emulate Theorem \ref{thm:measure0}, assuming
$\exists$-definability of $\Z$ in $\Q$ and building a c.e.\
set of products of primes so as to show that these are
the best possible bounds on the complexity of these real numbers.

\begin{thm}
\label{thm:rprime}
Assume that $\Z$ has an $\exists$-definition in $\Q$.
Then, given any two positive real numbers $u$ and $v$ with
$u+v < 1$, such that $u$ is computably approximable from
below and $v$ is $\HP$-computably approximable from below,
there exists a polynomial $f$ with $\alpha(f)=u$ and
$\gamma(f)=v$, hence with $\beta(f)=1-u-v$.
\end{thm}
\begin{pf}
We repeat the technique of Theorem \ref{thm:measure0}, by enumerating
the product $\Pi_{p\in I}~p$ of a finite set $I$ of primes into a c.e.\ set $D$
when we want the subring $\Z[I]$ to contain a solution to our polynomial $f$. 
This $f$ will be defined as $g^2+(XT-1)^2$ exactly as in that theorem,
using a polynomial $g(X,\Yvec)$ that defines $D$ in $\Q$ (which exists by the hypothesis
of $\exists$-definability of $\Z$ in $\Q$).  However, the enumeration of $D$ is
now more intricate: distinct elements of $D$ need no longer be relatively prime.

The enumeration of $D$ yields an enumeration of a c.e.\ set of nodes
$\sigma\in 2^{<\P}$:  those $\sigma$ such that
the products of the primes in $\sigma^{-1}(1)$ lies in $D$.
(Notice that a single element of $D$ may produce several such $\sigma$.
For example, if $35\in D$, then the strings $0011$, $1011$, $0111$,
and $1111$ all are enumerated.)
By the construction, then, these $\sigma$ will be precisely the nodes naming
the open set $\A(f)$, which will contain all subrings of $\Q$ extending any such $\sigma$.
At a stage $s$ in our construction, those $\sigma$ such that this product is
divisible by some $x\in D_s$ (that is, by some $x$ already enumerated into $D$)
will be said to be colored green at this stage.  (We think of them as
having  a ``green light'':  a solution to $f$ in the relevant subring is
already known.)  At stage $s$, the nodes colored red will be those nodes
$\tau$ such that no $\sigma\supseteq\tau$ is green.  Thus a node
may cease to be red at a particular stage, when it or a successor turns green;
if this happens, it will never again be red, although this node itself might
also never turn green.

By assumption there exsts a computable, strictly increasing
sequence $\la u_s\ra_{s\in\omega}$ of positive rational numbers
with $u=\lim_s u_s$.  Additionally, there exists a computable total
``chip'' function $c:\omega\to (0,1)\cap\Q$ such that
$$ \set{q\in\Q}{0<q<v} = \set{q\in\Q\cap (0,1)}{c^{-1}((0,q])\text{~is finite}},$$
so that the strict left Dedekind cut defined by $v$ is precisely
the set of rational numbers receiving only finitely many ``chips'' from $c$.
(\cite[Thm.\ IV.3.2]{S87} gives the essence of the construction of $c$.)
Notice also that with $v<1-u$, there will be infinitely many
$s$ with $c(s)<1-u<1-u_s$.  Indeed, by fixing a rational
number $q_0\in(v,1-u)$ and ignoring all stages $s$ with $c(s)>q_0$,
we may assume that $c(s)<1-u_s$ for every stage $s$.

At stage $0$, $D_0$ is empty.  At stage $s+1$, only finitely many
nodes can be minimal (under $\subseteq$) with the property of having
been green at stage $s$, since (by induction) $D_s$ was finite.  We fix the least level $l_s$
such that every minimal green node at stage $s$ lies at a level $\leq l_s$;
thus, at each level $\geq l_s$, every node must be either red or green at stage $s$.
(Below that level, a node may be neither color at stage $s$.)
We can list out the (finitely many) nodes that are minimal with the
property of having been red at stage $s$:  let these be $\rho_{0,s}\ldotc\rho_{j_s,s}$,
ordered by length so that $|\rho_{i,s}|\leq|\rho_{i+1,s}|$ and so that,
if these lengths are equal, then $\rho_{i,s}\prec\rho_{i+1,s}$ in the
lexicographic order $\prec$ on nodes.  We regard $\rho_{0,s}\ldotc\rho_{j_s,s}$
as a priority ordering of the minimal red nodes.

Recall that some rational $c(s+1)\in (0,1)$ received a chip
at this stage.  Find the greatest $k_s\leq j_s$ such that
$$ \sum_{i=0}^{k_s} 2^{-|\rho_{i,s}|} < c(s+1),$$
and for each of $\rho_{0,s}\ldotc\rho_{k_s,s}$, declare all of its
successors at level $l_s$ to be prioritized.  (This means that they
will all still be red at the end of this stage, and therefore so will
$\rho_{0,s}\ldotc\rho_{k_s,s}$.)  Let $\sigma_{0,s}\ldotc\sigma_{m_s,s}$
be the finitely many nodes of length $l$ that were red at stage $s$
but are not prioritized.

Our intention is to introduce a green node above each of these non-prioritized
nodes $\sigma_{i,s}$, so that they will no longer be red.
Therefore, for each $\sigma_{i,s}$ in turn, we enumerate into $D_{s+1}$
the product $x_{i,s}$ of a set of prime numbers such that:
\begin{itemize}
\item
whenever $\sigma_{i,s}(p)=1$, then $p$ divides $x_{i,s}$; and
\item
whenever $\sigma_{i,s}(p)=0$, then $p$ does not divide $x_{i,s}$; and
\item
$x_{i,s}$ has certain other prime factors $\notin\dom{\sigma}$, as defined below
(after Lemma \ref{lemma:red}).
\end{itemize}
The point of the first two rules is that now
$\sigma_{i,s}$ extends to some node that is colored green at stage $s+1$
(since $x_{i,s}\in D_{s+1}$).  However, we must ensure that no prioritized node
$\rho_{k,s}$ extends to a node that becomes green when $x_{i,s}$ enters $D$.
To understand why this is not immediate,
recall that if the number $35$ enters $D$, so as to
make the node $0011$ turn green, then the nodes $1011$, $0111$, and $1111$
will all also turn green, since they correspond to rings containing $\frac1{35}$.
However, Lemma \ref{lemma:red} shows
that these now-accidentally-green nodes cannot have been prioritized.
\begin{lemma}
\label{lemma:red}
For each $i\leq m_s$ and each $k\leq k_s$, some prime $q$ has $\rho_{k,s}(q)=0$
but $\sigma_{i,s}(q)=1$, so that $\frac1{x_{i,s}}\notin\Z[\P-\rho_{k,s}^{-1}(0)]$.
\end{lemma}
\begin{pf}
Let $\rho\subseteq\sigma_{i,s}$ be minimal such that $\rho$ is red
at stage $s$ (so, by our definition of $\sigma_{i,s}$ above,
$\rho = \rho_{j,s}$ for some $j>k_s$).  We must have
either $|\rho_{k,s}| < |\rho|$, or $\rho_{k,s}\prec \rho$.
Now if $\rho_{k,s}\prec \rho$, then the least
prime $q$ at which they differ has $\rho_{k,s}(q)=0$ and $\rho(q)=1$,
forcing $\sigma_{i,s}(q)=1$ since $\sigma_{i,s}$ restricts to $\rho$.
But if $|\rho_{k,s}| < |\rho|$, then $\rho\res |\rho_{k,s}|$ cannot be red as well,
by the minimality of $\rho$, and so some $q\in\rho_{k,s}^{-1}(0)$ must have
$\rho(q)=1$, for otherwise $\rho\res|\rho_{k,s}|$ would have been red
(as any green successor of $\rho\res|\rho_{k,s}|$ would have given rise
to a green successor of $\rho_{k,s}$).
\qed\end{pf}

By induction, we know that the measure $a_s$ of the set of all
paths in $2^{\P}$ that include a green node at stage $s$ lies
in $(u_s-\frac1{2^s},u_s)$, and we wish to make
$a_{s+1}\in (u_{s+1}-\frac1{2^{s+1}}, u_{s+1})$ as well.
(Recall that $u_s<u_{s+1}$.)
Now $a_s$ is precisely the measure of the set of nodes
at level $l_s$ that are green at stage $s$.  Meanwhile, the prioritized
nodes at level $l_s$ have total measure $<c(s+1)$, by our choice of $k_s$ above,
and these should stay red at stage $s+1$.  We arranged beforehand
that $u_{s+1}< 1-c(s+1)$, so that these requirements do not conflict.
The remaining nodes at level $l_s$ are precisely $\sigma_{0,s}\ldotc\sigma_{m_s,s}$.
Above we stated that each of these will contribute some $x_{i,s}$ to $D_{s+1}$.
By taking $x_{i,s}$ to have many prime factors $\notin\dom{\sigma_{i,s}}$,
we can make each $x_{i,s}$ contribute arbitrarily little measure to $a_{s+1}$,
so it is not difficult to ensure that $a_{s+1}<u_{s+1}$.  To make
$a_{s+1} > u_{s+1}-\frac1{2^{s+1}}$, we add a larger amount of measure
as needed, possibly enumerating several different numbers (but only finitely many)
into $D_{s+1}$ instead of just a single $x_{i,s}$.  For example, if $\sigma_{i,s}=0011$
(with $l_s=4$), then $5$ and $7$ must divide each $x_{i,s}$ and $2$ and $3$ must not;
by enumerating both $5\cdot 7\cdot 11\cdot 13$ and $5\cdot 7\cdot 11\cdot 17\cdot 19$ into $D_{s+1}$,
we can make the two extensions $0011\hat{~}11$ and $0011\hat{~}1011$ turn green.
(This would also make six other nodes, such as $1011\hat{~}11$, turn green,
if they were not green already.) 
These first two nodes together have measure $\frac5{16}\cdot\frac1{2^{l_s}}$, which
is five-sixteenths of the total measure available above $0011$.
How much else is added depends on whether $1011$, $0111$, and $1111$
were already green or not, but it is clear that we can compute this, and that
we could make any dyadic fraction of the total measure above the nodes $\sigma_{i,s}$
turn green.  So it is easy to make $a_{s+1}$ sit in the desired interval
$(u_{s+1}-\frac1{2^{s+1}},u_{s+1})$, effectively, and this completes the construction.

With $a_s\in (u_s-\frac1{2^s},u_s)$ for every $s$, it is clear that the resulting
polynomial $f$ has $\mu(\A(f))=\lim_s u_s = u$ as desired.  We also claim that
$\mu(\C(f))=v$, which will complete the proof.  In particular, whenever $q<v$,
we can produce a subset of $\C(f)$ of measure $\geq q$; whereas when $v<q$,
we will show that $\mu(\C(f))<q$ as well.

First suppose $q<v$, and fix any rational $q'\in (q,v)$.  Then there is some stage $s_0$
such that, for all $s\geq s_0$, we have $c(s)>q'$.   Then at stage $s_0$, among the 
minimal red nodes $\rho_{0,s_0}\ldotc\rho_{k_{s_0},s_0}$, the first $k$ (in this order)
will in fact be this highest-priority minimal red nodes remaining at the end of the construction,
where $k$ is maximal so that
$$ \sum_{i=0}^{k} 2^{-|\rho_{i,s_0}|} < q'.$$
Now $\rho_{k+1,s_0}$ may or may not remain red forever after.  If it does, then
we have a set of red nodes of total measure $\geq q'>q$, as required; so assume
that eventually a stage $s_1>s_0$ is reached at which some node above this $\rho_{k+1,s_0}$
is colored green.  Then at stage $s_1+1$, $\rho_{0,s_1+1}\ldotc\rho_{k,s_1+1}$ will
be the same as at stage $s_0$, but $\rho_{k+1,s_1+1}$ will be different:  either
it will have greater length than $\rho_{k+1,s_0}$, or it will be $\succ\rho_{k+1,s_0}$.
If it has the same length, then the same argument applies to this new $\rho_{k+1,s_1+1}$.
Since there are only finitely many nodes at each level, we either reach a node at this
level that stays red forever after, in which case again we have a set of red nodes of
sufficiently large measure; or else $\rho_{k+1,s}$ will eventually have greater length
than $\rho_{k+1,s_0}$.  This argument then continues until we reach a stage $s$ at which
$$2^{-|\rho_{k+1,s}|}<q'- \sum_{i=0}^{k} 2^{-|\rho_{i,s_0}|},$$
at which point this new $\rho_{k+1,s}$ will remain red forever.
The measure of the red set will become arbitrarily close to $q'$ via this process,
and hence must eventually be $>q$.  (With $q'<v$, it will eventually become $>q'$
as well, but this is irrelevant.)  To see why it must become arbitrarily close to $q'$,
notice that with $1-u_s>c(s)>q'$ at all subsequent stages, there will always be a supply
of red nodes of measure $>q'$, and the remainder of this measure will be partitioned
into smaller and smaller chunks as the length of the next minimal red node keeps increasing,
so that the measure of the permanently-minimal-red nodes cannot stay below $q'$
by any positive margin forever.

It remains to show that when $v<q$, we have $\mu(\C(f))<q$ as well.  Again it is useful
to fix some $q'$ between $q$ and $v$, now with $v<q'<q$.
Now there are infinitely many stages $s$ with $c(s)<q'$.  If the measure of $\C(f)$
were $>q'$, then evetually there would be a finite set of minimal red nodes, of total measure
$>q'$, all of which stayed red (and hence minimal) forever after.  But at some subsequent stage
$s$ we would have $c(s)<q'$, and at that stage the lowest-priority node in this finite set
would acquire a green node above it, so would not in fact have been permanently red.
With this contradiction, the proof is complete.
\qed\end{pf}

The remarks at the conclusion of Theorem \ref{thm:measure0} can be applied and expanded
here.  The first claim in this corollary follows from Theorem \ref{thm:measure0};
the second from Theorem \ref{thm:rprime}. 

\begin{cor}
\label{cor:computableboundary}
If the solution class $\A(f)$ of every polynomial $f\in\Z[\Xvec]$ has
computable measure, then there is no existential definition of $\Z$ in $\Q$.
Likewise, if the boundary class $\B(f)$ of every polynomial $f\in\Z[\Xvec]$ has
$\HP$-computable measure, then there is no existential definition of $\Z$ in $\Q$.
\qed\end{cor}

The following, while only a partial converse, serves to emphasize the
importance of the measures of boundary sets.

\begin{cor}
\label{cor:hardboundary}
If there exists a polynomial $f$ for which the measure $\beta(f)$ of $\B(f)$ is not
$\emptyset'$-computable -- or simply fails to be approximable from above --
then $HTP(\Q)$ is undecidable.
\end{cor}
\begin{pf}
If $HTP(\Q)$ is decidable, then the measures of both $\A(f)$ and $\C(f)$
are approximable from below, and therefore $\beta(f)=1-\alpha(f) - \gamma(f)$
is approximable from above.
\qed\end{pf}

\begin{cor}
\label{cor:Cfiniteunion}
Suppose that, for every polynomial $f\in\Z[X_0,X_1,\ldots]$, the set $\C(f)$ is an
effective union of basic open sets in $2^{\P}$.  (That is, suppose the red nodes
in $2^{\P}$ for $f$ always form a computably enumerable set.)  Then there is
no existential definition of $\Z$ in $\Q$.
\end{cor}
In particular, this corollary applies if, for each single $f$, the set of minimal red nodes for $f$
is a finite set.  The corollary would not require any method of determining the finite set
uniformly from the polynomial.  As of this writing, it is unknown whether there exists an $f$
for which the set of minimal red nodes is infinite (let alone not computably enumerable).
\begin{pf}
An effective union of basic open sets has as its measure a real number
approximable from below, and here this measure is $\gamma(f)$.
Since $\alpha(f)$ is always approximable from below, $\beta(f)=1-\alpha(f)-\gamma(f)$
would always be approximable from above, hence $\HP$-computable,
and we would then apply Theorem \ref{thm:rprime}.
\qed\end{pf}

\parbox{4.9in}{
{\sc
\noindent
Department of Mathematics \hfill \\
\hspace*{.1in}  Queens College -- C.U.N.Y. \hfill \\
\hspace*{.2in}  65-30 Kissena Blvd. \hfill \\
\hspace*{.3in}  Queens, New York  11367 U.S.A. \hfill \\
\ \\
PhD Programs in Mathematics \& Computer Science \hfill \\
\hspace*{.1in}  C.U.N.Y.\ Graduate Center\hfill \\
\hspace*{.2in}  365 Fifth Avenue \hfill \\
\medskip
\hspace*{.3in}  New York, New York  10016 U.S.A. \hfill}\\
\hspace*{.045in} {\it E-mail: }
\texttt{Russell.Miller\at {qc.cuny.edu} }\hfill \\
}


\begin{thebibliography}{99}

\comment{
\bibitem{C04}
S.\ Barry Cooper.
\newblock \emph{Computability Theory}.
\newblock Chapman \& Hall CRC, Boca Raton, FL, 2004.

\bibitem{C}
David A. Cox.
\newblock \emph{Primes of the Form $x^2+ny^2$:  Fermat, Class Field Theory and Complex Multiplication}.
\newblock John Wiley \& Sons, Inc., New York, 1989.
}

\bibitem{DPR61}
Martin Davis, Hilary Putnam, and Julia Robinson.
\newblock The decision problem for exponential diophantine equations.
\newblock {\em Annals of Mathematics}, 74(3): 425--436, 1961.

\bibitem{EMPS15}
Kirsten Eisentr\"ager, Russell Miller, Jennifer Park, and Alexandra Shlapentokh.
\newblock As easy as $\Q$:  Hilbert's Tenth Problem for subrings of the rationals.
\newblock \emph{Transactions of the American Mathematical Society}, 369(11): 8291--8315, 2017.

\bibitem{F57}
Richard M.\ Friedberg.
\newblock A criterion for completeness of degrees of unsolvability.
\newblock {\em Journal of Symbolic Logic}, 22: 159--160, 1957.

\bibitem{J81}
Carl G.\ Jockusch, Jr.
\newblock Degrees of generic sets.
\newblock {\em  Recursion theory: its generalisation and applications
(Proc. Logic Colloq., Univ. Leeds, Leeds, 1979)}, 
London Mathematical Society Lecture Note Series 45: 110--139, 1981.

\bibitem{JS83}
Carl G.\ Jockusch, Jr. and Richard A.\ Shore.
\newblock Pseudo jump operators I:  the r.e.\ case.
\newblock {\em Transactions of the AMS}, 275(2): 599--609, 1983.

\bibitem{JS84}
Carl G.\ Jockusch, Jr. and Richard A.\ Shore.
\newblock Pseudo-jump operators II:  transfinite iterations, hierarchies, and minimal covers.
\newblock {\em Journal of Symbolic Logic}, 49: 1205--1236, 1984.

\bibitem{K15}
Jochen Koenigsmann.
\newblock Defining {$\Z$} in {$\Q$}.
\newblock {\em Annals of Mathematics}, 183(1): 73--93, 2016.

\bibitem{KM19}
Kenneth Kramer and Russell Miller.
\newblock The Hilbert's-Tenth-Problem Operator.
\newblock \emph{Israel Journal of Mathematics}, 230(2): 693--713, 2019.

\bibitem{K81}
Stuart Kurtz.
\newblock \emph{Randomness and Genericity in the Degrees of Unsolvability.}
\newblock Ph.D.\ thesis, University of Illinois at Urbana-Champaign, 1981.


\bibitem{M70}
Yuri V. Matiyasevich.
\newblock The {D}iophantineness of enumerable sets.
\newblock {\em Dokl. Akad. Nauk SSSR}, 191: 279--282, 1970.


\bibitem{M08}
Russell Miller.
\newblock Computable fields and Galois theory.
\newblock Notices of the American Mathematical Society, 55(7): 798--807, 2008.

\bibitem{MCiE16}
Russell Miller.
\newblock Baire category theory and Hilbert's Tenth Problem inside $\Q$,
in \emph{Pursuit of the Universal:
12th Conference on Computability in Europe, CiE 2016},
eds.\ A.~Beckmann, L.~Bienvenu \& N.~Jonoska, Springer 
LNCS 9709: 343--352, 2016.

\bibitem{MIMS}
Russell Miller.
\newblock Measure theory and Hilbert's Tenth Problem inside $\mathbb Q$,
in \emph{Sets and Computations},
eds.\ S.D.\ Friedman, D.\ Raghavan, \& Y.\ Yang,
Institute for Mathematical Sciences, National University of Singapore,
Lecture Note Series 33: 253--269, 2017.

\bibitem{P03}
Bjorn Poonen.
\newblock Hilbert's Tenth Problem and Mazur's conjecture for large subrings of ${\Q}$.
\newblock {\em Journal of the AMS}, 16(4): 981--990, 2003.

\bibitem{P09}
Bjorn Poonen.
\newblock Characterizing integers among rational numbers with a universal-existential formula.
\newblock {\em American Journal of Mathematics}, 131(3): 675--682, 2009.

\bibitem{R49}
Julia Robinson.
\newblock Definability and decision problems in arithmetic.
\newblock {\em Journal of Symbolic Logic}, 14: 98--114, 1949.

\bibitem{R67}
Hartley Rogers.
\newblock{\em Theory of Recursive Functions and Effective Computability}.
McGraw-Hill, New York, 1967.

\bibitem{S73}
Jean-Pierre Serre.
\newblock \emph{A Course in Arithmetic}.
\newblock Graduate Texts in Mathematics, vol.\ 7 (Berlin: Springer, 1973).

\bibitem{S07}
Alexandra Shlapentokh.
\newblock{\em Hilbert's Tenth Problem: Diophantine Classes and Extensions to Global Fields}.
Cambridge U.P., Cambridge, UK, 2007.

\bibitem{S87}
Robert I. Soare.
\newblock {\em Recursively Enumerable Sets and Degrees}.
\newblock Perspectives in Mathematical Logic. Springer-Verlag, Berlin, 1987.



\end{thebibliography}
\end{document}